\documentclass[11pt, oneside]{article}   	
\usepackage{bussproofs}
\usepackage{geometry}                		
\usepackage{graphicx}
\usepackage{tikz}
\usepackage{pifont}
\usepackage{linguex}
\usepackage{natbib}

\usepackage{amssymb}
\usepackage{amsmath}

\usepackage{hyperref}
\hypersetup{colorlinks=true,allcolors=blue}

\newcommand{\F}{\mathsf{F}}
\newcommand{\s}{\mathsf{S}}

\newcommand{\PE}[1]{{\color{black} #1}}

\title{Truth and Falsity in Buridan's Bridge}

\author{Paul \'Egr\'e}

\date{}							

\begin{document}
\maketitle

\sloppy

\begin{abstract}  

\noindent This paper revisits Buridan's Bridge paradox (\emph{Sophismata}, chapter 8, Sophism 17), itself close kin to the Liar paradox, a version of which also appears in Bradwardine's \emph{Insolubilia}. Prompted by the occurrence of the paradox in Cervantes's \emph{Don Quixote}, I discuss and compare four distinct solutions to the problem, namely \PE{Bradwardine's ``just false'' conception, Buridan's ``contingently true/false'' theory, Cervantes's ``both true and false'' view, and then Jacquette's ``neither true simpliciter nor false simpliciter'' account. All have in common to accept that the Bridge expresses a truth-apt proposition, but only the latter three endorse the transparency of truth}. Against some previous commentaries I first show that Buridan's solution is fully compliant with an account of the paradox within classical logic. I then show that Cervantes's insights, as well as Jacquette's treatment, are both supportive of a dialetheist account, and Jacquette's in particular of the strict-tolerant account of truth. I defend dialetheist intuitions (whether in LP or ST guise) against two objections: one concerning the future, the other concerning the alleged simplicity of the Bridge compared to the Liar.

\end{abstract}


\section{Introduction}

Faced with the antinomy of the Liar (``this sentence is not true''), one option is to judge that upon further analysis the sentence should be considered neither true nor false. Arguably, the sentence fails to express any genuine proposition. Various analogies can help here. Maybe the sentence is like ``Colorless green ideas sleep furiously", it is grammatical but we cannot assign it any clear meaning. Or maybe the sentence is like ``The present King of France is bald'', it lacks a presupposition and cannot be evaluated as either true or false for that matter. Or maybe it is like ``she was born in Mexico'': if you don't know to whom the pronoun ``she'' refers, then the meaning of the sentence is incomplete, and the sentence cannot be evaluated as true or false, it needs further input from the context.\footnote{See in particular \cite{bochvar1937, fraassen1968presupposition,glanzberg2004contextual} for ways in which these analogies (with meaninglessness, presupposition failure, context-dependence) can be fleshed out and articulated.}

But is it so clear that the Liar sentence is neither true nor false? Dialetheists like Graham Priest think otherwise, and view the sentence as both true and false (\citealt{priest1979logic,priest2006incontradiction}). One way in which dialetheism can be defended is by considering an analogy between the Liar and vague sentences (\citealt{mcgee1990truth, cobreros2015vagueness,egre2019respects, egre2021half}). Unlike absurd sentences, presupposition failures, or sentences with unresolved indexical pronouns, the Liar may be viewed as a borderline case of application of the vague predicate ``true''. Vague sentences are not like absurd sentences, nor like cases of presupposition failure (\citealt{Zehr2014:PHD, spector2016multivalent}). And while they can be viewed as indexical sentences of a specific kind, they admit judgments of truth and falsity in borderline cases, implying that indexical elements in them are resolved. 
According to various analyses of vagueness at least, vague sentences in borderline cases are sentences that are true in some sense, but also false in some sense, in that regard they do not lack truth values (viz. \citealt{hyde1997heaps, cobreros2012tcs}).

Intuitions about the Liar and related paradoxes are hard to trust and to make stable, however, and the justification of dialetheism on grounds of vagueness is controversial (see \citealt{priest2019respectfully}). Hence, one may wonder if one can find independent evidence for the judgment that the problematic sentence is both true and false. Toward that goal, in this paper I propose to revisit a variant of the Liar paradox commonly known as Buridan's Bridge.  
The Bridge is close kin to the Liar, and indeed Buridan presents it and discusses it in Chapter 8 of his \emph{Sophismata} (under Sophism 17), where he also examines the Liar (under Sophism 11). The puzzle, in fact, already appears in Bradwardine's treatise on \emph{Insolubilia} (8.8.1), where Bradwardine handles it in the same way in which he handles the Liar. In recent years, however, Buridan's Bridge appears to have received less attention than the Liar, at least compared to the medieval and Renaissance period for which there is evidence of its popularity (see \citealt{ashworth1972treatment, ashworth1976bridge, jones1986liar}), but also compared to the vast exegesis on Buridan's treatment of the Liar (see \citealt{prior1962problems}, \citealt{read2002liar}, and further references below). Among the papers that focus specifically on the sophism of the Bridge and on Buridan's solution, the most recent sources I could find include \citealt{hughes1982buridan}, \citealt{jacquette1991buridan}, and \citealt{ulatowski2003buridan}.\footnote{Mention of the puzzle can be also found in \citealt{prior1962problems}. Prior in that paper does not propose an analysis of that specific puzzle, however. Instead he focuses on earlier puzzles in the \emph{Sophismata}, in particular the Liar. More recent sources citing the paradox include \citealt{read1995thinking} (about Cervantes), and \citealt[36]{ohrstrom2007temporal} (about Buridan).} None of them, however, offers a formal semantic treatment of the puzzle. 

One motivation for this paper is to fill this gap, but more fundamentally, the main reason to look at the Bridge is that it is a case for which one can find assessments that the sentence is both true and false, though coming from sources that are not usually associated with the dialetheist tradition. In particular, a version of Buridan's Bridge can be found in Cervantes' \emph{Don Quixote} (part II, chap. 51).\footnote{See \citet[283]{prior1962problems}, who writes: ``this is a puzzle of some literary interest, since there is a very similar one in \emph{Don Quixote}''. Prior cites \citet[Exercice 15.10]{church1996introduction} regarding the connection to Cervantes. Church states a version of the puzzle very close to Cervantes's original, to which he gives credit, but he does not refer to Buridan. As I will argue below, the treatment sketched by Church happens to agree with a central element in the diagnosis proposed by Buridan himself, although Church does not appear to be aware of Buridan's account. I myself discovered Church's treatment only after formalizing Buridan's account. No reference to either Church or Prior appears in the texts by Hughes, Jacquette, or Ulatowski. More generally, the literature on the Bridge paradox is very scattered.} In the novel, Sancho Panza is asked to issue a verdict on the case, and his verdict is that there is as much truth in the problematic sentence as there is falsity in it. As I will argue, I think Cervantes's discussion should be taken seriously, and that it lends support to dialetheist intuitions.\footnote{\PE{\citet[p. 185, fn. 5]{priest2006incontradiction}, in his chapter on norms and the philosophy of law, does not mention the Bridge. However, he briefly mentions the related case of Protagoras and Euathlus (see below, fn. \ref{fn:sources}) as a case of \textit{prima facie} dialetheia, but he does not elaborate on it.
}
} 
The second source is the analysis of the paradox proposed by Dale Jacquette. While Jacquette remains classically-minded throughout his paper, his informal remarks bring him so close to dialetheism that, combined with Cervantes's naive observations, both accounts should be considered as valuable data points toward showing that dialetheist intuitions can arise outside of a dialetheist framework when it comes to semantic paradoxes.

Before showing this, however, I propose to examine Buridan's paradox closely. As a first step, I will argue that the simplicity of Buridan's own solution was missed by \PE{some} previous commentators. Jacquette, particularly, writes: 
``If Buridan's remarks are taken at face value, they are difficult to reconcile with contemporary ways of thinking about these matters'' \citep[460]{jacquette1991buridan}. \PE{However, Jacquette used \cite{buridanscott}'s edition and translation of the \emph{Sophismata}, which turns out to be unreliable in several places. In one passage, Scott's translation includes a wrong negation, which misled Jacquette into thinking that Buridan analyzes conditionals in a nonclassical fashion.}\footnote{\PE{I am indebted to an anonymous referee for bringing this fact to my notice. On pg. 458 of his paper, Jacquette reproduces Scott's translation, which says: ``I say that Plato does not speak truly, since Socrates has \emph{not} fulfilled the condition" (emphasis mine, based on the reviewer's observation). In the editions of \cite{hughes1982buridan,klimaburidan,buridanpironnet}, the problematic negation is missing, viz. \cite{klimaburidan}: ``I say that Plato did not say something true, since Socrates fulfilled the
condition''.}
} 
\PE{Relying on later editions \citep{hughes1982buridan, klimaburidan,buridanpironnet}, one can formally reconstruct Buridan's arguments to see that his analysis is fully compliant with classical logic.}

Buridan's own solution to the Bridge, precisely because of its simplicity, may suggest that there is no need for a nonclassical treatment of the puzzle, let alone a dialetheist approach. However, I would like to show that as much in the Bridge paradox as in the Liar paradox, a dialetheist treatment is natural if we insist on making all premises jointly true together. To establish this, I shall rely on some of Jacquette's own suggestions. Ironically, Jacquette dismisses Cervantes's treatment of the Bridge.\footnote{See \citet[457]{jacquette1991buridan}: ``Sancho gives an alternative philosophically less satisfactory resolution in this strange but no less edifying History, which it were now tedious to relate'', and his remark ``Cervantes is evidently a better novelist than logician."} Yet he suggests that a completely classical analysis of the paradox can be given, despite considering that the Bridge sentence
\begin{quote}
is neither true \emph{simpliciter}, nor false \emph{simpliciter} ({which is not to say that it is neither true nor false}) for demonstrably, the sentence is true if and only if it is false.
\end{quote}
I will argue that Jacquette's informal remarks can be reconciled with Cervantes' own diagnosis, and that they can also be used to support to the variety of dialetheist analysis exemplified in the strict-tolerant analysis of the Liar (\citealt{cobreros2013reaching,ripley2012conservatively,ripley2013paradoxes}). The fact that Jacquette sees his own diagnosis as compatible with classical logic may not be an accident. 
It exemplifies the way in which strict-tolerant logic can be made to agree with classical logic, by emphasizing the fact that the predicate ``true'' can be taken in a weak or a strong sense.

The main steps in the paper are as follows: I first review the presentations of the paradox given by Bradwardine, by Buridan, and by Cervantes, to highlight their common structure, and to give a brief overview of their respective solutions. I then formalize the paradox in first-order logic. The formalization is used to successively clarify Buridan's account (section \ref{sec:buridan}), Cervantes's account (section \ref{sec:cervantes}), and finally Jacquette's account (section \ref{sec:jacquette}). \PE{All three have in common to endorse the transparency of truth, namely the equivalence between  ``$A$'' and ``$\langle A\rangle$ is true'', a principle denied in Bradwardine's antecedent approach}. I then examine a revenge objection against the LP and ST accounts, concerning the semantics of the future (section \ref{sec:future}). A second objection concerns the prima facie superiority of a classical treatment of the Bridge in comparison to its alternatives (section \ref{sec:liar}). I argue that this verdict lacks generality, and conclude with broader lessons regarding the comparison between the classical, the LP, and the ST account of the Bridge.

\section{\PE{The} Bridge}\label{sec:paradox}

\subsection{The paradox}

The Bridge paradox is stated as Sophism 17 in Chapter 8 of Buridan's \emph{Sophismata} (written in the first half of the fourteenth century).
In Buridan's version, Plato is master of a bridge and will not let anyone pass without his permission:\footnote{This is \cite{hughes1982buridan}'s translation, slightly modified, based on the comparison with Klima's translation \citep{klimaburidan} and the Latin text edited by Pironnet \citep{buridanpironnet}.}

\begin{quote}

Then Socrates arrives on the scene and pleads urgently with Plato to let him cross. 
Then Plato, irritated, makes an oath and says: ``Surely Socrates, if in the first proposition you utter you speak truly, I will let you cross; but be sure that if you speak falsely, I will throw you in the water''.

\end{quote}

\noindent Socrates's response consists in the utterance: ``you will throw me in the water'' (\emph{proiicies me in aquam}). Buridan's presentation of the paradox (or sophisma) then goes as follows:

\begin{quote}
If it is said that he will throw Socrates in the water, it goes against his promise, for then Socrates spoke truly, and Plato should let him cross. If it is said that he lets him cross, then it goes again against the oath and promise, for then Socrates spoke falsely, in which case Plato should have thrown him in the water.

\end{quote}

\PE{
Buridan likely did not invent the puzzle, even though his version is special and his analysis original. However, the paradox has earlier roots in Ancient Philosophy, in what is sometimes referred to as the Crocodile dilemma, and in Protagoras' paradox.\footnote{\label{fn:sources}See \citet[161]{hughes1982buridan}. The Crocodile dilemma: ``a crocodile stole a woman's baby and promised to return it to her if she told him truly whether he would it eat it or not. The woman then replied, ``you are going to eat it''." The mother argues that whether what she says is true or false, the Crocodile should keep his word and return the baby; the crocodile argues symmetrically.  Protagoras's paradox: Protagoras sues his student Euathlus, who had promised to pay his fee when winning his first case, but chose never to practice. Protagoras argues that if he loses, Euathlus will have to pay him (because he lost), and if he wins, Euathlus will also have to pay him (to keep his word); Euathlus reasons that if he loses he won't have to pay (short of winning his first case), and similarly if he wins (because then Protagoras should pay). Hughes does not give the actual sources, but the Crocodile dilemma is attributed to Chrysippus by Lucian of Samosata in his \emph{Vitarum Auctio} (see \citealt{biard1993sophismata}, as well as \citealt[part III, commentary to chap. 8, p. 159]{murray1847compendium}). The dilemma of Protagoras and Euathlus, also cited in \citet[part III, chap. 8, p. 156]{murray1847compendium}, is mentioned in Diogenes La\"ertius \emph{Vitae philosophorum}, IX, 56 and in Aulus Gellius, \emph{Noctes atticae} V, 10. The case is also presented by Sextus Empiricus in \emph{Adversus Rhetores} (€96-100), except that Sextus mentions Korax (a Sicilian orator linked to Gorgias) instead of Protagoras, and an unnamed disciple instead of Euathlus. In Sextus's version, the judges eventually kick Korax and his student out of the court, on the grounds that their arguments have equal strength. Walker in his commentary of Murray's compendium does not cite Buridan, but notes the connection of the ancient paradoxes with Cervantes: ``if I mistake not, a similarly ludicrous instance may be found in Don Quixote''. 
}} \PE{The puzzle also appears in other treatises on insolubles by several contemporaries of Buridan's, in particular in Bradwardine's \emph{Insolubilia} \citep[p. 124-25; 134-35]{bradwardineinsol}, generally assumed to precede Buridan's treatise (ca. 1325) and to have influenced him. Bradwardine's treatise even contains two variants of the problem. In the first (8.1), the main hypothesis is that ``all and only those who do not utter a falsehood will receive a penny'', and Socrates declares ``Socrates will not receive a penny''. The second variant (8.8.1) corresponds to the bridge case proper, with main hypothesis that ``all and only those speaking the truth will cross the bridge'', and Socrates declares ``Socrates will not cross the bridge''.\footnote{\PE{Pseudo-Heytesbury's \emph{Insolubilia Padua} also presents the penny version \citep[294]{pironnet2008william}, namely ``Socrates non habebit denarium". See also Wyclif's \emph{Summa insolubilium} (ca. 1360), Book III, chapter 13 [261], which refers to it as ``casu[s] de pertransitione pontis'' \citep{wyclifsumma}. I am indebted to an anonymous referee for these references.}}}

The version given by Cervantes in \emph{Don Quixote}, while coming three centuries later, is very similar in structure (\citealt{quixote}, part II, chap. 51), even though Cervantes does not credit any external source or authority.
\footnote{\cite{jones1986liar} conjectures that Cervantes might have relied on treatises by the Spanish philosophers Domingo de Soto (\emph{Introductiones Dialectice}, 1529) and Juan de Celaya (\emph{Insolubilia et Obligationes}, 1517), as well as on Paul of Venice's \emph{Logica Magna}, first published in 1499. Surprisingly, Jones does not cite Buridan \PE{or earlier sources}. On the analyses of the Bridge given by Paul of Venice and by Domingo de Soto, see \cite{ashworth1976bridge}.} In the story, Sancho Panza is a newly appointed governor of the Island of Barataria, and he is challenged with a series of problematic legal cases. In Cervantes' presentation, four judges stand before a bridge and have to ensure that:

\begin{quote}
If anyone crosses by this bridge from one side to the other he shall declare on oath where he is going to and with what object; and if he swears truly, he shall be allowed to pass, but if falsely, he shall be put to death for it by hanging on the gallows erected there, without any remission.
\end{quote}

A man comes up and declares by oath that ``he was going to die upon that gallows that stood there, and for nothing else'' (\emph{iba a morir en aquella horca que all\'i estaba, y no a otra cosa}). Cervantes's presentation of the antinomy faced by the judges parallels the one given by Buridan:

\begin{quote}

If we let this man pass free he has sworn falsely, and by the law he ought to die; but if we hang him, as he swore he was going to die on that gallows, and therefore swore the truth, by the same law he ought to go free.

\end{quote} 

Admittedly, the \PE{different} 
versions mentioned here differ in some details. \PE{Bradwardine's version involves Socrates making a third-person assertion about himself}. In Buridan's version, Socrates's utterance is reported in direct speech, it states a future tense sentence about a second-person action.  Cervantes's version is a third-person indirect report, in the past tense, of a first-person utterance made in present tense but concerning the future. I will set aside those grammatical differences here, to focus mostly on the common logical structure behind \PE{those} formulations.\footnote{\cite{jacquette1991buridan} considers that Cervantes's version makes the paradox easy to solve, because Cervantes's version refers to the bystander's first-person \emph{intention} to cross the bridge. I believe Jacquette's assessment is based in part on an artefact of the English translation (namely the translation of ``y no a otra cosa'' by ``and that was all his business''), which can be ignored.}

\subsection{Bradwardine's, Buridan's, and Cervantes's solutions}

\PE{
In the sequel, I will be mostly interested in the comparison between Buridan's solution and Cervantes's approach to the puzzle, as both of which subscribe to the transparency of truth (the intuitive idea that a sentence and the assertion of its truth can be exchanged in non-intensional contexts).\footnote{The principle is also called the Intersubstitutivity Principle in \citealt{field2008saving}, and Na\"ivet\'e in \citealt{rossi2019unified}.} To appreciate their originality, however, it is useful to consider Bradwardine's treatment first, which does not endorse transparency. 

For Bradwardine, Socrates's utterance ``Socrates will not cross the bridge'' is just false. The reason, according to Bradwardine, is \citep[8.8.1, p. 135]{bradwardineinsol}:

\begin{quote}

If he will cross, then he utters a falsehood, so he will not cross, so it must be said that Socrates will not cross. On the contrary: then [Socrates's utterance] is true, that is, Socrates spoke the truth. It must be said that it does not follow.

\end{quote}

\noindent That is, Bradwardine accepts the argument that Socrates's crossing would make his utterance false, 
implying he will not cross (given the main hypothesis, namely the decree). But he rejects the inference from ``Socrates will not cross the bridge'' to the conclusion that the sentence is true. The reason advanced by Bradwardine, developed after the passage cited above, is that the decree makes ``crossing the bridge'' and ``speaking the truth'' equivalent as a matter of fact. So for Bradwardine, all that the sentence \emph{signifies} is that Socrates will not cross, which is equivalent in context to the fact that the sentence is false.\footnote{Bradwardine's argument rests on a principle explained in Chapter 6 of his \emph{Insolubles}, namely ``every proposition signifies or means as a matter of fact or absolutely respectively everything which follows from it as a matter of fact or absolutely''. Note that I leave aside a formalization of Bradwardine's theory in what follows, to focus on the other three accounts, which subscribe to transparency.} 
His analysis therefore makes two very specific assumptions: it rejects the ascent from ``$A$'' to ``$\langle A\rangle$ is true'' (one side to transparency), and moreover it takes the truth of the decree for granted.}

Buridan's analysis is very different from Bradwardine's, since for Buridan Socrates' utterance is not false \emph{simpliciter}, but contingently true or false.
According to Buridan, Socrates's utterance is either true or false, but it is not known which until Plato acts one way or the other. As Buridan also puts it: Socrates's utterance fails to be determinately true or to be determinately false. Buridan's main reason is that Socrates's utterance concerns a future contingent, which it is in Plato's power to make true or false. Buridan, furthermore, infers that Plato's oath must be false, and that Plato will have to violate his promise, which is another important difference compared to Bradwardine's treatment in particular.

In Cervantes's narrative, on the other hand, Sancho is requested to give a verdict on whether the bystander spoke truly or falsely. Sancho's reaction is in two steps. While both of his judgments may appear to be merely poking fun of scholastic reasoning, they are worth taking more seriously. Sancho, after all, impersonates a character whose intelligence is meant to be spontaneous and vivid, despite his lack of instruction (or even \emph{because} of his lack of instruction -- as also exemplified in Descartes' democratic conception of common sense a few years later). 

Sancho's first proposal regarding the case is to say that:

\begin{quote}

of this man they should let pass the part that has sworn truly, and hang the part that has lied; and in this way the conditions of the passage will be fully complied with.

\end{quote}

When Sancho's interlocutor points out to him that one can not divide a man without causing him to die, Sancho comes up with this second verdict:

\begin{quote}

either I'm a numskull or else there is the same reason for this passenger dying as for his living and passing over the bridge; for \emph{if the truth saves him the falsehood equally condemns him}. And that being the case it is my opinion you should say to the gentlemen who sent you to me that \emph{as the arguments for condemning him and for absolving him are exactly balanced}, they should let him pass freely, as it is always more praiseworthy to do good than to do evil; this I would give signed with my name if I knew how to sign. (emphasis mine)

\end{quote}

As those citations indicate, logically speaking Sancho considers that the Bridge sentence is both true and false, or at any rate, that it is partly true, and partly false. Pragmatically, Sancho's verdict is to let the man pass free, but this tolerant choice is driven by ethical considerations that are external to the argument. If the man over the bridge turned out to be an awful murderer about to kill innocent victims, Sancho could faultlessly choose to hang the man by giving salience to the fact that the utterance is false, on the same consequentialist grounds (namely to minimize evil and maximize goodness).

In summary, Buridan treats the Bridge sentence as (just) true or (just) false, and views Plato's decree as (just) false. Cervantes on the other hand takes the truth of the decree for granted, and ends up with the conclusion that the Bridge sentence is both true and false. \PE{Both of their approaches differ from Bradwardine's, who also takes the truth of the decree for granted, but argues that the Bridge sentence is just false, thereby rejecting the transparency of truth.}

\section{The Bridge, formally}\label{sec:formalism}

\subsection{Set-up}

Before engaging further into semantic analysis, let us handle the paradox formally. To do so I shall use a first-order language (without identity), enriched with a future tense operator. Let $Sxy$ represent the predicate ``$x$ utters $y$'', ``$Tx$'' represent ``$x$ is true'', ``$Px$'' represent ``$x$ is hanged/thrown in the water/does not cross the bridge'', and let $\F$ represent the operator ``it will be the case that''. To represent the rules proclaimed in either story (the decree/the oath), without loss of generality we may rely on the version proposed by Cervantes and use the following formulae T and F:

\ex.[(T)] \a.[] {Who speaks truly will not be hanged/thrown in the water/cross the bridge}.\label{true}
\b.[] $\forall x\forall y(Sxy \to (Ty \rightarrow \neg \mathsf{F}Px))$

\ex.[(F)] \a.[] {Who speaks falsely will be hanged/thrown in the water/not cross the bridge}.\label{false}
\b.[] $\forall x \forall y (Sxy \to (\neg Ty \rightarrow \F Px))$

\noindent The conditionals expressed in the oath are here treated as truth-functional material conditionals, in agreement with Buridan's own comments on the way in which the conditional can be interpreted in this case.\footnote{Buridan contrasts two analyses of conditionals in his discussion of the second question raised by the puzzle, one we may call strict (a conditional cannot be true if it is possible for the antecedent to be true and the consequent false), and one we can call material, based on the fact that it obeys the principle of conjunctive sufficiency (a conditional is true when antecedent and consequent are both true). Buridan argues that the latter suffices for the argument. \PE{On the tradition to treat promissive conditionals as material conditionals until the Renaissance period, see in particular \cite{ashworth1972strict}.}} \PE{The treatment of the conditionals as truth-functional can furthermore justify the idea that Plato's oath can be meaningfully called true or false. Intuitively, an oath of promise is either \emph{fulfilled} or \emph{unfulfilled}, but it does not appear to be truth-apt. However, Socrates' utterance about the future is truth-apt, since it makes a forecasting use of the future: either what it says will happen as it says, or it won't. Likewise, consider Plato's oath in Buridan's version: either Plato will let Socrates cross as a matter of fact, or he will throw him in the water. Hence, as a truth-functional compound, Plato's promissive conditional is itself truth-apt, since the antecedent and the consequent of each of its conjuncts are truth-apt. Henceforth, we make the same assumption about the universally quantified conditionals (T) and (F), namely that they can be verified or falsified.\footnote{\PE{In other words, I handle ``will'' throughout in a predictive sense, and not in the deontic sense of ``ought'' in which it is sometimes interpreted. Throughout this paper I therefore talk of the decree as being true or false, instead of fulfilled vs. unfulfilled. One may, in principle, introduce a predicate of fulfillment distinct from the truth predicate, in order to give a more elaborate account. For more on the distinction between truth and fulfillment, see \cite{bonevac1990paradoxes}. For a defense of the view that legal statements are evaluable as true and false, see \citet[p.~186]{priest2006incontradiction}.}}}

Let us use $\phi$ as an abbreviation for the conjunction of F and T, corresponding to Plato's oath/the river lord's decree. Let us use $a$ to denote Socrates/the bystander, and let us name $b$ the Bridge sentence, which can be represented as follows:

\ex.[($b$)] \a.[] {I will be hanged/thrown in the water.}
\b.[] $\F Pa$

To model truth, we shall use the following two principles of transparency, where $\langle A\rangle $ is a name of the sentence $A$:

\ex.[(Tr)] Transparency of truth
\[
\begin{tabular}{ccc}
{$T\langle  A \rangle$} & & {$A$}\\
\cline{1-1}  \cline{3-3}
$A$ & & $T\langle A\rangle$\\
\end{tabular}
\]

Finally, I assume that these principles are preserved under negation, namely:

\ex.[($\neg$Tr)] Transparency of truth under negation:
\[
\begin{tabular}{ccc}
{$\neg T\langle A \rangle$} & & {$\neg A$}\\
\cline{1-1}  \cline{3-3}
$\neg A$ & & $\neg T\langle A \rangle $\\
\end{tabular}
\]

Other choices would be appropriate. One may use the biconditional schema $T\langle A\rangle \leftrightarrow A$. However, it will be simpler in what follows to represent Transparency in rule form.

\subsection{Derivation of the antinomy}

Several derivations of the antinomy may be given. In this section, I start with a derivation in natural deduction style which is faithful to the informal reasoning followed by Buridan. I also present it this way because it helps illuminate Buridan's philosophical analysis of the paradox. The bolded formulae in the tree are meant to highlight the part of the argument that is made explicit in the citations given above.

\ex.\label{tree:buridan} Reasoning from future facts

\hspace{-1cm} \scalebox{.7}{
\parbox{1cm}{
\begin{prooftree}

\AxiomC{$\F Pa \vee \neg \F Pa$}

\AxiomC{$\phi$}
\LeftLabel{\scriptsize $\wedge$E to T}
\UnaryInfC{$\forall x \forall y(Sxy \to (Ty \to \neg \F Px))$}
\LeftLabel{\scriptsize $\forall$E}
\UnaryInfC{$Sab \to (Tb \to \neg \F Pa)$}

\AxiomC{$Sab$}
\LeftLabel{\scriptsize $\to$E}
\BinaryInfC{$\boldsymbol{Tb} \to \neg \F \boldsymbol{Pa}$}

\AxiomC{$[\F \boldsymbol{Pa}]$}
\LeftLabel{\scriptsize Tr}
\UnaryInfC{$\boldsymbol{Tb}$}
\LeftLabel{\scriptsize $\to$E}
\BinaryInfC{$\neg \F \boldsymbol{Pa}$}
\LeftLabel{\scriptsize $\wedge$I}
\UnaryInfC{$ \F Pa \wedge \neg \F Pa$}

\AxiomC{$\phi$}
\LeftLabel{\scriptsize $\wedge$E to F}
\UnaryInfC{$\forall x \forall y(Sxy \to (\neg Ty \to \F Px))$}
\LeftLabel{\scriptsize $\forall$E}
\UnaryInfC{$Sab \to (\neg Tb \to \F Pa)$}

\AxiomC{$Sab$}
\LeftLabel{\scriptsize $\to$E}
\BinaryInfC{$\neg \boldsymbol{Tb} \to \F \boldsymbol{Pa}$}

\AxiomC{$[\neg \F \boldsymbol{Pa}]$}
\LeftLabel{\scriptsize $\neg$Tr}
\UnaryInfC{$\neg \boldsymbol{Tb}$}
\LeftLabel{\scriptsize $\to$E}
\BinaryInfC{$ \F \boldsymbol{Pa}$}
\LeftLabel{\scriptsize $\wedge$I}
\UnaryInfC{$ \F Pa \wedge \neg \F Pa$}

\LeftLabel{\scriptsize $\vee$E}
\TrinaryInfC{$\F Pa \wedge \neg \F Pa$}

\end{prooftree}
}}\bigskip

Plato's oath (hypothesis $\phi$) and the fact that Socrates utters the Bridge sentence ($Sab$) are taken for granted. In both Buridan's version and in Cervantes's version the main assumptions concern what will happen in the future, namely the factual hypotheses $\F Pa$ and $\neg \F Pa$ used to reason by cases. From both assumptions, using Transparency (modulo negation), a contradiction follows, namely that Socrates will be thrown in the water, and that Socrates will not be thrown in the water. 

An equivalent derivation can be produced, this time by reasoning by case, metalinguistically, upon the truth of the Bridge sentence. This is the same derivation, except that $\F Pa$ and $Tb$ are swapped to begin with. This alternative derivation will prove useful further below.

\ex.\label{tree:truth} Reasoning from the truth of the Bridge sentence

\noindent \scalebox{.7}{
\parbox{1cm}{
\begin{prooftree}

\AxiomC{$Tb \vee \neg Tb$}

\AxiomC{$\phi$}
\LeftLabel{\scriptsize $\wedge$E to T}
\UnaryInfC{$\forall x \forall y(Sxy \to (Ty \to \neg \F Px))$}
\LeftLabel{\scriptsize $\forall$E}
\UnaryInfC{$Sab \to (Tb \to \neg \F Pa)$}

\AxiomC{$Sab$}
\LeftLabel{\scriptsize $\to$E}
\BinaryInfC{$Tb \to \neg \F Pa$}

\AxiomC{$[Tb]$}
\LeftLabel{\scriptsize $\to$E}
\BinaryInfC{$\neg \F Pa$}
\LeftLabel{\scriptsize $\neg$Tr}
\UnaryInfC{$ \neg Tb$}
\LeftLabel{\scriptsize $\wedge$I}
\UnaryInfC{$ Tb \wedge \neg Tb$}

\AxiomC{$\phi$}
\LeftLabel{\scriptsize $\wedge$E to F}
\UnaryInfC{$\forall x \forall y(Sxy \to (\neg Ty \to \F Px))$}
\LeftLabel{\scriptsize $\forall$E}
\UnaryInfC{$Sab \to (\neg Tb \to \F Pa)$}

\AxiomC{$Sab$}
\LeftLabel{\scriptsize $\to$E}
\BinaryInfC{$\neg Tb \to \F Pa$}

\AxiomC{$[\neg Tb]$}
\LeftLabel{\scriptsize $\to$E}
\BinaryInfC{$ \F Pa$}
\LeftLabel{\scriptsize Tr}
\UnaryInfC{$Tb$}
\LeftLabel{\scriptsize $\wedge$I}
\UnaryInfC{$Tb \wedge \neg Tb$}

\LeftLabel{\scriptsize $\vee$E}
\TrinaryInfC{$Tb \wedge \neg Tb$}

\end{prooftree}
}}


\subsection{Noteworthy features}

Two remarks can be made about those derivations, in comparison to other paradoxes based on self-reference in particular. Firstly, a noteworthy feature is that $b$ does not say of itself directly that it is not true. Principles T and F produce the same effect but in an indirect way, since under the assumption $Sab$, $\F Pa$ and $\neg Tb$ are interderivable. 

Secondly, it is not clear whether the future tense operator $\F$ plays any substantive role. As can be checked, the operator may be erased from the statement of T, F, and $b$, then in the trees, without impugning the soundness of the derivation. Does this mean that consideration of time and tense should be deemed irrelevant?

Prima facie it does indeed, and this fact may be used to vindicate Jacquette's analysis, according to which considerations of time and tense are a red herring. There is a subtlety, however, when we consider Sancho's verdict. Suppose the man is let free to pass over the bridge and not hanged. Does that not make the Bridge sentence $b$ true after all? But if so, how can the paradox be avoided? This may be seen as a form of revenge. As I will argue below, Buridan's solution takes it into account, but it does not appear to be seen by Cervantes in his discussion.

\section{Buridan's solution: classicalism}\label{sec:buridan}

Famously, Buridan considers three questions in relation to the paradox. The first is whether Socrates' utterance is true or false. The second is whether Plato's oath is true or false. The third concerns what Plato must do to keep his oath. The three answers have been rather sternly criticized, in particular by Jacquette and by Ulatowski.\footnote{\citet[459]{jacquette1991buridan} writes: ``The analysis jumbles several things together. Again there are apparent semantic errors from the standpoint of modern logic". \citet[87]{ulatowski2003buridan} hammers the point: ``...Buridan's solutions to the bridge paradox are inadequate: the logico-semantic issue jumbles several things together...". \PE{However, both Jacquette and Ulatowski actually rely on Scott's edition and translation of Buridan. Hughes' edition explicitly differs from Scott's in a number of places, and Hughes correctly explains how Plato's conditional must be false (\citealt[160-161]{hughes1982buridan}), yet without formalizing the argument.}} I shall argue that their criticisms are unfair: neither of them formalizes the puzzle, and they miss the classicality of Buridan's solution.

Buridan's answer to the first question is that Socrates's answer cannot be \emph{known} to be true or false before observing Plato's action. For him, it is in Plato's power to make it true or false. Importantly, this does not mean that Buridan considers that Socrates's utterance is neither true nor false. What Buridan writes is that Socrates' utterance is not determinately true or determinately false at the time of the utterance, but Buridan is explicit that Socrates ``has uttered a proposition that must be true or false''. Hence, this assessment is compatible with Buridan holding that the law of excluded middle applies in the case of $b$, namely $\F Pa \vee \neg \F Pa$ must hold true for him.

This point is crucial, and in fact can be used to explain his answer to the second question. His answer to the second question is that Plato's oath is not true (\emph{Plato non dixit verum}), and therefore must be false (\emph{sequitur ipsam esse falsam}).\footnote{\PE{Wyclif's treatment of the Bridge is much closer to Buridan's than it is to Bradwardine's in that regard. See \cite[261-262]{wyclifsumma}, where Wyclif writes that even if a tyrant makes a universal order as expressed in the bridge decree, ``non ex hoc quod ipse sic ordinat sequitur universalem istam esse veram'' (from the fact that he orders things that way, it does not follow that this universal statement is true).}}
 The simplest way to vindicate this reasoning is to see that given that Buridan accepts that one of $\F Pa$ and $\neg \F Pa$ must be true, then given the noncontroversial truth of $Sab$, and Buridan's endorsement of the instances of Transparency in the derivation, the only assumption left to challenge is $\phi$, namely Plato's oath. In other words, Buridan appears to simply use a classical \emph{reductio} argument.

\ex. Buridan's reductio against $\phi$

\hspace{-1cm} \scalebox{.7}{
\parbox{1cm}{
\begin{prooftree}

\AxiomC{$\F Pa \vee \neg \F Pa$}

\AxiomC{$[\phi]$}
\LeftLabel{\scriptsize $\wedge$E to T}
\UnaryInfC{$\forall x \forall y(Sxy \to (Ty \to \neg \F Px))$}
\LeftLabel{\scriptsize $\forall$E}
\UnaryInfC{$Sab \to (Tb \to \neg \F Pa)$}

\AxiomC{$Sab$}
\LeftLabel{\scriptsize $\to$E}
\BinaryInfC{$\boldsymbol{Tb} \to \neg \F \boldsymbol{Pa}$}

\AxiomC{$[\F \boldsymbol{Pa}]$}
\LeftLabel{\scriptsize Tr}
\UnaryInfC{$\boldsymbol{Tb}$}
\LeftLabel{\scriptsize $\to$E}
\BinaryInfC{$\neg \F \boldsymbol{Pa}$}
\LeftLabel{\scriptsize $\wedge$I}
\UnaryInfC{$ \F Pa \wedge \neg \F Pa$}

\AxiomC{$[\phi]$}
\LeftLabel{\scriptsize $\wedge$E to F}
\UnaryInfC{$\forall x \forall y(Sxy \to (\neg Ty \to \F Px))$}
\LeftLabel{\scriptsize $\forall$E}
\UnaryInfC{$Sab \to (\neg Tb \to \F Pa)$}

\AxiomC{$Sab$}
\LeftLabel{\scriptsize $\to$E}
\BinaryInfC{$\neg \boldsymbol{Tb} \to \F \boldsymbol{Pa}$}

\AxiomC{$[\neg \F \boldsymbol{Pa}]$}
\LeftLabel{\scriptsize $\neg$Tr}
\UnaryInfC{$\neg \boldsymbol{Tb}$}
\LeftLabel{\scriptsize $\to$E}
\BinaryInfC{$ \F \boldsymbol{Pa}$}
\LeftLabel{\scriptsize $\wedge$I}
\UnaryInfC{$ \F Pa \wedge \neg \F Pa$}

\LeftLabel{\scriptsize $\vee$E}
\TrinaryInfC{$\F Pa \wedge \neg \F Pa$}
\LeftLabel{\scriptsize EFQ}
\UnaryInfC{$\bot$}
\LeftLabel{\scriptsize $\to$I}
\UnaryInfC{$\neg \phi$}

\end{prooftree}
}}\bigskip

Jacquette writes: ``If Buridan's remarks are taken at face value, they are difficult to reconcile with contemporary ways of thinking about these matters''. This is, however, an uncharitable reading of Buridan. A careful examination of Buridan's argument reveals that his solution is perfectly classical indeed. As a further confirmation, \cite{church1996introduction} presents a treatment of Cervantes's version of the paradox in propositional logic whose conclusion is exactly the same, namely that the law must be false, given the other assumptions. While Church may not even have been aware of Buridan's account, his own formal treatment vindicates exactly Buridan's second conclusion.\footnote{Church's treatment in Exercise 15.10 runs as follows: $P$ represents the fact that the man crosses the bridge, $R$ the fact that the man's utterance is true, $Q$ the fact that the man is hanged, and $S$ the fact that the law is obeyed. His assumptions are that $R\leftrightarrow P \wedge Q$, $P$, and $S\to (Q\leftrightarrow (P \wedge \neg R))$. Church's exercise is to show that $\neg S$ follows from those assumptions, by modus ponens and given the tautology $(R\leftrightarrow P \wedge Q) \to (P\to ((S\to (Q\leftrightarrow P \wedge \neg R))\to \neg S))$.}

To the third question, finally, Buridan answers that Plato does not have to hold his promise. The reason is first stated in relation to the second question, namely ``because of Socrates' proposition Plato's promise has reference to
itself whence it follows that it is false
". In plain terms, Plato's promise is self-refuting due to Socrates' utterance. Buridan continues and writes (\citealt{klimaburidan}'s translation, slightly modified): 

\begin{quote}

he does not have to keep his promise, nor should he promise
anything in this way, but with an exception that excludes the case that Socrates
utters a proposition that has reference to the promise such that it thence
follows that what is promised cannot be made true.


\end{quote}

Buridan's analysis stops here. To clarify his point, we can observe that if $\phi$ is false, then by classical reasoning, its negation must be true. The negation of $\phi$ is equivalent to the disjunction of the negations of principles T and F, namely:

\ex. $\exists x \exists y (Sxy \wedge Ty \wedge \F Px) \vee \exists x \exists y (Sxy \wedge Fy \wedge \neg \F Px)$

That is, either someone utters a true sentence and will be thrown in the water, or someone utters a false sentence and will not be thrown in the water. 

Now, let us suppose that Socrates is eventually thrown in the water by Plato. Then what Socrates said turns out to have been true, which instantiates the first disjunct. And suppose that Socrates is not thrown in the water. Then what Socrates said turns out to have been false, which instantiates the second disjunct.\footnote{This point is acknowledged by Jacquette, see \cite[462]{jacquette1991buridan}. Although Jacquette sees this, he does not appear to see it as the conclusion of a classical reductio argument.} This confirms that at the moment Socrates made his utterance, one of $Tb$ and $\neg Tb$ was true indeed, though it was not determinate which. But Plato's acting one way or the other prevents fulfilling his promise.

\section{Cervantes's solution: dialetheism}\label{sec:cervantes}

Instead of rejecting the decree as false, like Buridan, in Cervantes's novel Sancho takes a naive approach, namely presumes that the decree is true, and takes the conclusion of the reasoning in \ref{tree:buridan} at face value. Since Sancho ends up with the observation that there are equal reasons to hang the man and to let him pass, one way to regiment this view is to adopt an LP-style approach (\citealt{priest1979logic}), namely to declare $b$ to be both true and false. 


In a three-valued setting, a straightforward model of Sancho's reasoning is therefore to let $\F Pa$ take the value $\frac{1}{2}$, and to assign $Sab$ the value $1$, since the fact that the bystander uttered the Bridge sentence is just true. An issue concerns the choice of a semantics for the conditional. One option is to adopt the Strong Kleene conditional as in LP. In that case all steps in \ref{tree:buridan} receive a value at least $\frac{1}{2}$. With the LP conditional, conditional elimination is no longer valid, however, given the definition of LP-validity as the preservation of non-zero values. 

To make all steps valid, what is needed is a conditional supporting modus ponens. This could be a modal conditional (\citealt{priest2006incontradiction}), or a different three-valued conditional, such as the Cooper conditional (\citealt{cooper1968propositional}), which takes the value of the consequent when the antecedent has a value 1 or $\frac{1}{2}$, and the value $\frac{1}{2}$ when the antecedent is false.\footnote{See \citealt{egre2021finettian} for a presentation. Whether the Cooper conditional is suitable to deal with other self-referential paradoxes beside the Bridge is not obvious, however, but I leave this issue aside here.} For conjunction, the most natural choice is to use the min-conjunction, and to evaluate the universal quantifier as the generalization of that operator. For negation, we assume the strong Kleene negation.

Under those assumptions, namely assuming a modus-ponens-supporting conditional combined with LP-validity, all steps in \ref{tree:buridan} come out as valid and all can be sound. The important point is that the decree $\phi$ gets the value $\frac{1}{2}$, given that the instance $Sab \to (Tb \to \neg \F Pa)$ gets the value $\frac{1}{2}$. This means that $\phi$ too is both true and false, like its negation.

Practically, if the man passes over the bridge free, as Sancho recommends, this is justified by the fact that the sentence is true indeed, given that it is true and false. But if the man were hanged at the gallows, this would be justified by the fact that the sentence was also false, because true and false.\footnote{\citet[149]{read1995thinking} comments on Sancho's decision to let the man pass free as follows: ``this is \emph{in effect} to set the law aside, and declare it inoperable in this case. In other words, the law should have been more carefully framed in the first place, and was unsatisfactory'' (my emphasis). That is, Read looks at the decision from the viewpoint of classical logic. But note that this does not represent Sancho's point of view. Sancho himself does not claim that the law is inoperable.}

Despite its simplicity, the dialetheist solution may appear repugnant, since one may object that $\F Pa$ and its negation ought to take a determinate value 1 or 0. After all, if Socrates passes over the bridge eventually, how can it be the case, retrospectively,  that Socrates passes and does not pass over the bridge? I return to that issue below.

\section{Jacquette's solution: strict vs tolerant}\label{sec:jacquette}

A close kin to the dialetheist solution just presented is Jacquette's treatment of the puzzle. According to him (\citealt[466]{jacquette1991buridan}, emphasis in the original):

\begin{quote}

Plato can either permit Socrates to pass or have him seized and thrown into the river \emph{without violating his conditional vow}. The reason is that Plato swore only what he would do \emph{if} Socrates' proposition was true (\emph{simpliciter}), and what he would do \emph{if} Socrates' proposition was false (\emph{simpliciter}), but he left open what he would do or would not do if Socrates' proposition was \emph{not} true \emph{simpliciter} or false \emph{simpliciter}. 
\end{quote}

Jacquette's main argument toward that conclusion is that Socrates's utterance is neither true simpliciter nor false simpliciter, since it is provably true if and only if false given Plato's decree. What follows is a formal construal of Jacquette's argument, for which we introduce the operator $\s$, whereby $\s Tx$ represents that $x$ is true \emph{simpliciter}. We call ``the $\s$-rule'' the rule used by Jacquette, according to which if a sentence is true iff false, then it is neither true simpliciter nor false simpliciter.

\ex. Jacquette's argument

\scalebox{.75}{
\parbox{1cm}{
\begin{prooftree}

\AxiomC{$\phi$}
\LeftLabel{\scriptsize $\wedge$E}
\UnaryInfC{$\forall x \forall y(Sxy \to (Ty \to \neg \F Px)$}
\LeftLabel{\scriptsize $\forall$E}
\UnaryInfC{$Sab \to (Tb \to \neg \F Pa)$}

\AxiomC{$Sab$}
\LeftLabel{\scriptsize $\to$E}
\BinaryInfC{$Tb \to \neg \F Pa$}

\AxiomC{$Tb$}
\LeftLabel{\scriptsize $\to$E}
\BinaryInfC{$\neg \F Pa$}
\LeftLabel{\scriptsize $\neg$Tr}
\UnaryInfC{$ \neg Tb$}
\LeftLabel{\scriptsize $\to$I}
\UnaryInfC{$ Tb \to \neg Tb$}

\AxiomC{$\phi$}
\LeftLabel{\scriptsize $\wedge$E}
\UnaryInfC{$\forall x \forall y(Sxy \to (\neg Ty \to \F Px)$}
\LeftLabel{\scriptsize $\forall$E}
\UnaryInfC{$Sab \to (\neg Tb \to \F Pa)$}

\AxiomC{$Sab$}
\LeftLabel{\scriptsize $\to$E}
\BinaryInfC{$\neg Tb \to \F Pa$}

\AxiomC{$\neg Tb$}
\LeftLabel{\scriptsize $\to$ E}
\BinaryInfC{$ \F Pa$}
\LeftLabel{\scriptsize Tr}
\UnaryInfC{$Tb$}
\LeftLabel{\scriptsize $\to$I}
\UnaryInfC{$\neg Tb \to Tb$}

\LeftLabel{\scriptsize $\leftrightarrow$I}
\BinaryInfC{$Tb \leftrightarrow \neg Tb$}
\LeftLabel{\scriptsize $\s$-rule}
\UnaryInfC{$\neg \s Tb \wedge \neg \s \neg Tb$}

\end{prooftree}
}}

Jacquette's argument here is arguably compatible with a dialetheist treatment, since dialetheism too \PE{can distinguish} between sentences that are just true and sentences that are true and false.\footnote{\PE{This claim is controversial, as highlighted to me by David Ripley, depending on how ``just true'' is formalized. In \cite{beall2009spandrels}'s extension of LP with transparent truth, $T\langle A\rangle$ and $T\langle A\rangle \wedge \neg T\neg \langle A \rangle$ are equivalent. If the latter is taken to paraphrase ``just true'', then ``true'' and ``just true'' collapse.}} One way to implement the argument semantically is to handle the $\s$ operator as a Bochvar-type meta-assertion operator (\citealt{bochvar1937,beaver2001}). That is, given a three-valued model $v$, let $v(\s A)=1$ if $v(A)=1$, and $v(\s A)=0$ otherwise. This suffices to make the $\s$-rule LP-valid.

According to Jacquette, what Plato really intends by his decree is the conjunction we can call $\phi^{\star}$ of:

\ex.[(T$^\star$)] \a.[] {Who speaks truly \emph{simpliciter} will not be hanged/thrown in the water}.\label{truestar}
\b.[] $\forall x\forall y(Sxy \to (\s Ty \rightarrow \neg \mathsf{F}Px))$

\ex.[(F$^\star$)] \a.[] {Who speaks falsely \emph{simpliciter} will be hanged/thrown in the water}.\label{falsestar}
\b.[] $\forall x \forall y (Sxy \to (\s \neg Ty \rightarrow \F Px))$

\noindent As Jacquette admits, however (\citealt[468]{jacquette1991buridan}), this diagnosis rests on pragmatic reasoning, for it requires an inference about Plato's communicative intentions:

\begin{quote}
The solution requires that we assume Plato would agree that his intent in communicating the vow was to declare what he would do if Socrates' first utterance were true \emph{simpliciter}, and what he would alternatively do if Socrates' first utterance were false \emph{simpliciter}.
\end{quote}

In relation to this remark, we may also look at Jacquette's analysis within the strict-tolerant framework of \cite{cobreros2013reaching}, in which truth comes in two levels, a tolerant and a strict level. In ST, a sentence $A$ is tolerantly true provided $v(A)\geq \frac{1}{2}$, and a sentence is strictly true provided $v(A)=1$. Using a Kleene-style conditional, we have to remember that $A\to B$ is tolerantly true provided if $A$ is strictly true, then $B$ is tolerantly true (viz. \citealt{cobreros2015vagueness}). Viewed that way, the conditionals (T$^{\star}$) and (F$^{\star}$) both come close to tolerant readings of rules (T) and (F), since the tolerant reading of the conditional requires the truth predicate to be interpreted strictly in each occurrence.\footnote{Moreover, (T$^{\star}$) and (F$^{\star}$) come out equivalent to tolerant readings of (T) and (F) in three-valued logic if we assume that for $Sxy$ and $\F Pa$ the distinction between tolerant and strict collapses.} It is therefore natural to equate truth \emph{simpliciter} with strict truth.

An alternative presentation of Jacquette's argument that does not use the $\s$-operator explicitly is obtained by extending \ref{tree:truth} into a rejection of the Law of Excluded Middle for $b$:

\ex.\label{tree:st} Reductio against LEM

\noindent \scalebox{.7}{
\parbox{1cm}{
\begin{prooftree}

\AxiomC{$[Tb \vee \neg Tb]$}

\AxiomC{$\phi$}
\LeftLabel{\scriptsize $\wedge$E to T}
\UnaryInfC{$\forall x \forall y(Sxy \to (Ty \to \neg \F Px))$}
\LeftLabel{\scriptsize $\forall$E}
\UnaryInfC{$Sab \to (Tb \to \neg \F Pa)$}

\AxiomC{$Sab$}
\LeftLabel{\scriptsize $\to$E}
\BinaryInfC{$Tb \to \neg \F Pa$}

\AxiomC{$[Tb]$}
\LeftLabel{\scriptsize $\to$E}
\BinaryInfC{$\neg \F Pa$}
\LeftLabel{\scriptsize $\neg$Tr}
\UnaryInfC{$ \neg Tb$}
\LeftLabel{\scriptsize $\wedge$I}
\UnaryInfC{$ Tb \wedge \neg Tb$}

\AxiomC{$\phi$}
\LeftLabel{\scriptsize $\wedge$E to F}
\UnaryInfC{$\forall x \forall y(Sxy \to (\neg Ty \to \F Px))$}
\LeftLabel{\scriptsize $\forall$E}
\UnaryInfC{$Sab \to (\neg Tb \to \F Pa)$}

\AxiomC{$Sab$}
\LeftLabel{\scriptsize $\to$E}
\BinaryInfC{$\neg Tb \to \F Pa$}

\AxiomC{$[\neg Tb]$}
\LeftLabel{\scriptsize $\to$E}
\BinaryInfC{$ \F Pa$}
\LeftLabel{\scriptsize Tr}
\UnaryInfC{$Tb$}
\LeftLabel{\scriptsize $\wedge$I}
\UnaryInfC{$Tb \wedge \neg Tb$}

\LeftLabel{\scriptsize $\vee$E}
\TrinaryInfC{$Tb \wedge \neg Tb$}
\LeftLabel{\scriptsize EFQ}
\UnaryInfC{$\bot$}
\LeftLabel{\scriptsize $\to$I}
\UnaryInfC{$\neg (Tb \vee \neg Tb)$}

\end{prooftree}
}}

Using the strong Kleene conditional, this derivation is no longer LP-valid, due to the failure of modus ponens in LP, and also due to the reliance on the Ex Falso Quodlibet rule. {Each step is  $ss$-valid (preserving strict truth), \PE{except for the final discharge.\footnote{Throughout we read a natural deduction proof of $A$ with set of open assumptions $\Gamma$ as a proof of $\Gamma, \Delta\models A$, for any $\Delta$. The metainference from $\Gamma, A\models^{ss} \bot$ to  $\Gamma\models^{ss} \neg A$ is not $ss$-valid. (Thanks to D. Ripley for a related discussion).
}} {Every step is $st$-valid, however, meaning that if each premise is assumed to hold strictly, the next conclusion follows tolerantly. In particular, $Tb \vee \neg Tb$ cannot be assumed to hold strictly, assuming $\phi$ to hold strictly. However, the negation of $Tb \vee \neg Tb$ can hold tolerantly}.

What about the decree, in this approach? {In ST, the reductio rule remains in place, and so $\neg \phi$ too can be inferred from the previous derivation}. This means that the negation of $\phi$ holds tolerantly, under the assumption that $\phi$ holds strictly. This is compatible with the negation of $\phi$ holding strictly, but also with $\phi$ holding tolerantly. In other words, the ST approach need not decide whether the negation of $\phi$ holds strictly, or whether it only holds tolerantly. As I will emphasize in the conclusion, I think this flexibility is a strength of the ST framework.

Jacquette in his account leaves open two interpretations: either the Bridge sentence fails to express any proposition, or  

\begin{quote}
if not all genuine propositions need to be true \emph{simpliciter} or false \emph{simpliciter}, but at least some genuine propositions are true if and only if they are false, then Socrates may have expressed a genuine proposition by this utterance
\end{quote}

\noindent Clearly, however, Jacquette's account leans toward the latter interpretation. While Jacquette refrains from committing to any non-classical theory in his account, he therefore admits a weaker notion of truth than the notion of truth \emph{simpliciter}.

I conclude that Jacquette's account is compatible with an LP-style version of dialetheism as much as with a close kin like the strict-tolerant account of truth, since both of these recognize that not all genuine propositions need to be true or false \emph{simpliciter}. This is not to say that Jacquette would endorse the strict-tolerant account, given that the latter involves accepting that logical consequence can fail to be transitive (viz. \citealt{ripley2013paradoxes}). But arguably, Jacquette's attachment to classical inferences, combined with the observation that truth \emph{simpliciter} is pragmatically conveyed in the Bridge, can be used to favor the latter.

\section{The future}\label{sec:future}

We may now compare the classical solution with the LP and ST solutions. Prima facie, the classical solution appears to be free of a revenge problem that confronts dialetheism. For suppose that Socrates is set free to cross the bridge, as Sancho's clemency recommends. Then Socrates's utterance turns out false, since it is not the case that he is thrown in the water by Plato. But if so, does this not make $\neg \F Pa$ true \emph{simpliciter} as a result? The same reasoning obviously applies in case Socrates is thrown in the water: then it appears $\F Pa$, and thus $b$, ends up true \emph{simpliciter}.

This is a problem for a non-classical account of the paradox along the lines of LP or ST. Is there a way around it? To solve it, one option is to distinguish the fact that in a three-valued setting, $\F Pa$ takes the value $\frac{1}{2}$ at time $t$, from the fact that $Pa$ takes a classical value 1 or 0 at a later time. Indeed, we run into trouble only if we assume a plain Priorean semantics for the future:

\ex. $v(t,\F A)=1$ iff $\exists t'>t: v(t',A)=1$

To avoid this problem, another option is to assume a Thomason-style supervaluationist semantics for the future, by relativizing the evaluation of the $\F$ operator to histories (\citealt{thomason1970indeterminist,cobreros2016supervaluationism}):

\ex. \a. $v(t,\F A)=1$ if {for all histories $h$} through $t$, there is $t'>t$ such that $v(h,t', A)=1$
\b. $v(t, \F A)=0$ if {for all histories $h$} through $t$, there is $t'>t$ such that $v(h,t',A)=0$
\c. $v(t, \F A)=\frac{1}{2}$ otherwise.

We may use these truth conditions for the cases in which $A$ is itself a sentence containing no $\F$ operator or $T$ predicate. For complex sentences built out of those sentences, again the strong Kleene operators may be assumed. Assume that at $t$, two different histories $h$ and $h'$ are such that at a later time $t'$, $v(h,t',Pa)=1$ and $v(h',t',Pa)=0$. This reflects the indeterminacy of Plato's action relative to time $t$. This implies that $v(\PE{t},\F Pa \vee \neg \F Pa)=\frac{1}{2}$. This agrees with the assignment of the value $\frac{1}{2}$ to the disjunction $Tb \vee \neg Tb$ at time $t$. 
It is therefore compatible with a dialetheist account of the Bridge sentence $b$ and of Plato's promise $\phi$ to assume that at a time posterior to Socrates's utterance, Plato will take a determinate action, assigning $Pa$ a classical truth value. 

As a result, the main difference here with Buridan's treatment is that Buridan regards the disjunction $\F Pa \vee \neg \F Pa$ as true \emph{simpliciter}, whereas the present account handles this disjunction as itself taking an intermediate truth value, namely as only tolerantly true. This is another step away from Buridan's approach, in particular from his endorsement of the principle of bivalence. As we saw earlier, the assumption that $\F Pa \vee \neg \F Pa$ is true \emph{simpliciter} was a pivotal element in order to conclude that Plato's promise was false. Whether in LP or ST, however, the Law of Excluded Middle remains valid, but it is not true strictly.

\section{The Liar within the Bridge}\label{sec:liar}

Despite the coherence of the LP and ST approaches, a further objection that may be made to both is that the classical solution is just simpler. Indeed, from Buridan's solution, one may be tempted to conclude that the Bridge is only a ``so-called paradox'', to use Quine's terminology (\citealt{quine1953so}), hence not a real antinomy, since it simply rests on a specious premise. But are things so simple?

Let us consider the question from Buridan's perspective first. While commentators rightly classify the group of sophisms among which the Bridge appears as ``pragmatic paradoxes" in Buridan's Chapter 8 (namely paradoxes dealing with promises, wishes, intentions, or knowledge),\footnote{See \citealt{prior1962problems, burge1978buridan, hughes1982buridan}.} it is not immediately obvious whether Buridan viewed the Bridge and the Liar as fundamentally different species. What is common to both treatments is that Buridan deals with the Liar paradox using a similar \emph{reductio} strategy to the one he uses in the Bridge. On the Liar, however, Buridan's reductio faces a difficulty which forces him to give up the transparency of truth, just like we saw Bradwardine do in the Bridge.
\footnote{See \cite{read2002liar} for a detailed account of Bradwardine's treatment of the Liar, and its likely influence on Buridan.}  

Thus, in Sophism 11, Buridan first concludes that the Liar sentence (``what I am saying is false'') is {false}, again because assuming its truth leads to a contradiction, namely to the conclusion that the sentence is both true and false. The crux of his discussion is then to block the inference from the Liar's falsity to the conclusion that the Liar is true (the semantic ascent from ``$A$'' to ``$\langle A \rangle$ is true'', which is half of transparency). Here is Buridan's objection to that move (trad. Hughes, modified):

\begin{quote}

It is claimed that if the proposition is false it follows that it is true. I reject that inference. You may try to defend it on the ground that if the proposition is false then the facts are as it says they are. That I admit, as far as its formal meaning is concerned. But this is not enough to make it true, due to the reflection upon itself that it comports. And in fact it is \emph{not} true, since the facts are not as the conclusion implied in it and the case would say they are.


\end{quote}  

While the similarity to Bradwardine's treatment is clear, this passage is notoriously difficult (see \citealt{prior1962problems,herzberger1973dimensions,hughes1982buridan, read2002liar, klima2008logic, novaes2009lessons,perini2011buridan} for critical examinations), and it would lie beyond the scope of this paper to engage into a precise discussion of Buridan's conception of truth.
\footnote{
As in the case of the Bridge, I believe that the principle of charity should be rigorously applied. Presently, however, I find harder to accept Buridan's (and Bradwardine's) denial that the Liar's falsity should entail is truth. Despite that, two abstract elements of convergence may actually be underlined concerning Buridan's treatment of the Liar and the Bridge, and the strict-tolerant account in particular: both accounts accept that those sentences are truth-apt, and both accept that truth is not a simple notion. In Buridan, it involves distinct dimensions, see \cite{herzberger1973dimensions}; in ST, it comes in a weak vs a strong form.} What we see, however, is that for Buridan, the Liar is a self-refuting sentence (the supposition of its truth leads to contradiction), but that like Bradwardine, Buridan needs to limit transparency in order to block the paradox. In the case of the Bridge, Buridan does not need that move because he can reject the truth of the decree instead, an option Bradwardine does not entertain.

However, one may use the problematic status of the Liar to argue that the Bridge is too easily solved by simply rejecting Plato's oath as being false.
For instance, one may ask what Plato should do in the Bridge situation if Socrates had uttered the Liar sentence ``what I am saying is false'' instead of the Bridge sentence ``you will throw me in the water''. Given his analysis of the Liar, Buridan would conclude that Socrates should be thrown in the water in this case, since the Liar is false \emph{simpliciter} according to him. But what if we doubt his conclusion that the Liar is false without being thereby true? What if, like Sancho with the Bridge, we find compelling the argument that the falsity of the Liar should entail its truth?

Granting that the Liar is true if and only if false, then we are back to a situation in which Socrates uttered a sentence that is neither true \emph{simpliciter}, nor false \emph{simpliciter}. The difference with the Bridge in its original formulation is that the truth or falsity of Socrates's modified utterance (i.e. the Liar) does not depend, this time, on what action Plato will choose to perform. 

This variation on the puzzle can be used to buttress the intuition that Plato's oath is not fundamentally the culprit in Buridan's Bridge. The problem, in agreement with Jacquette's analysis, concerns the existence of sentences that are neither true \emph{simpliciter} nor false \emph{simpliciter}. 



\section{Concluding remarks}

We have compared distinct solutions to the Bridge paradox, all of which accept 
the fact that the Bridge sentence expresses a genuine proposition, evaluable as true or false. First Bradwardine's solution, which concludes that the Bridge sentence is just false, but at the cost of rejecting transparency; then,
Buridan's solution, which concludes that Plato's oath must be false, but handles the Bridge as contingently true or false
; then Cervantes' dialetheist solution, which admits that the Bridge sentence is true and false, and gives the same status to Plato's oath; finally, the strict-tolerant solution, which agrees with the classical solution on the fact that Plato's oath is not strictly true, but which also agrees with the dialetheist account on the fact that both the Bridge and Plato's oath can be handled as tolerantly true. 

Let us highlight a few lessons from this comparison. The first is that Cervantes's discussion as well as Jacquette's discussion may be viewed as ``natural cases'' in support of a dialetheist account of semantic paradoxes. The existence of such ``natural cases'' is not just anecdotal, since Jacquette's account in particular gives reasons to view the Bridge sentence as neither true \emph{simpliciter} nor false \emph{simpliciter}, while granting that the sentence is also true and false. \PE{Further support may actually be drawn from medieval logic in this regard. 
Very interestingly, a reviewer points out that Richard Kilvington, in his \emph{Sophismata} (written in the first quarter of the fourteenth century), puts forward a view of insolubles that anticipates both Jacquette's remarks and the strict-tolerant account of semantic paradoxes. Kilvington, in his Sophisma 48, writes in particular: ``I say, then, that no insoluble is absolutely true or absolutely false; instead each is true in a certain respect and false in a certain respect'' (\citealt[142]{kilvingtonsophismata}). Kilvington's discussion there deals explicitly with the Liar and with a passage of Aristotle's \emph{Sophistical Refutations} in which Aristotle asks whether someone swearing only to perjure himself swears truly or not (180b-181a).\footnote{For a discussion of this passage, in particular in relation to Bradwardine's treatment of the Liar, see \cite{dutilh2008insolubilia}.} The crux of Kilvington's discussion is the observation that whereas ``true absolutely'' and ``false absolutely'' are contraries, the terms ``true in a certain respect'' and ``false in a certain respect" are not. While a careful examination of Kilvington's conception lies beyond this paper, it sets the strict-tolerant distinction concerning ``true'' and ``false'' in an even broader perspective and tradition.}

The second lesson from our analysis is exegetical, namely that Buridan's own account can be fully regimented within classical logic. Church independently reached the same verdict as Buridan's in his presentation of Cervantes's version, and Hughes acknowledged the consistency of Buridan's analysis, but the fact that Buridan's second conclusion can be vindicated by classical logic was missed by several commentators, despite the common element between Buridan's treatment of the Liar and his treatment of the Bridge. {\PE{In that respect, our account confirms a remark made by \citet[p.~201]{perini2011buridan}, which is that ``nowhere in his \textit{Sophisms} does Buridan say, or even suggest, that a sentence may lack truth-value''. In the case of the Bridge, this is definitely so, and this undermines \citet[p.~91]{ulatowski2003buridan}'s claim that ``Buridan implicitly accepts a three-valued logic when he claims Socrates' proposition is a future contingent''. As we saw, Buridan's analysis complies entirely with bivalence, whereas the introduction of a third value in our account is first and foremost motivated by the intuition that the Bridge and the decree could both be true together.}}

The third and main point is that the simplicity of Buridan's solution to the Bridge paradox does not suffice to recommend his solution. Buridan basically viewed the Liar and the Bridge as cases for which a reductio argument should lead us to the conclusion that a particular sentence (the Liar itself, Plato's oath) must be false \emph{simpliciter}. But beside the fact that Buridan's verdict on the Liar is moot, I have argued that by replacing the Bridge sentence by the Liar sentence, the same paradoxical conclusion is obtained. The fundamental problem in the Bridge, as much as in the Liar, concerns the existence of sentences that are true without being true \emph{simpliciter} (in ST terms, sentences that are tolerantly true, without being strictly true).

Finally, in this paper I have not undertaken a deeper comparison between the LP account of the Bridge and the ST account. Admittedly, the two frameworks have much in common, and ST has even been described by some as ``LP in sheep's clothing''.\footnote{See \citealt{barrio2015logics,dicher2019st,cobreros2020inferences,priest2021substructural} for comparisons and discussions of that claim.} 
Despite their proximity, I think a strength of the ST approach is the fact that it can preserve Buridan's diagnosis that the oath cannot be true in the strict sense, but compatibly with the dialetheist approach, precisely based on the recognition that truth can be interpreted either weakly or strongly. In that regard, and even setting aside the controversial account of logical consequence it encapsulates, the strict-tolerant approach does represent a synthesis between a standard classical account, and a standard dialetheist account.

\section*{Acknowledgments}

{\footnotesize I wish to thank Jean-Pascal Anfray, Pablo Cobreros, Dimitri El Murr, Julie Goncharova, Louis Guichard, James Hampton, Klaus von Heusinger, Andrea Iacona, Benjamin Icard, Hitoshi Omori, Paloma P\'erez Ilzarbe, Simone Picenni, Graham Priest, Lorenzo Rossi, David Ripley, Johannes Stern, and Denis Thouard for very helpful exchanges, and audiences in G\"ottingen (AIL1), Bristol, Torino, Pamplona, and Paris. Special thanks go to two anonymous referees and to the editors of this special issue for comments. I am particularly grateful to Paloma P\'erez Ilzarbe and to David Ripley for their detailed feedback and generous comments on specific aspects of the paper, to Jean-Pascal Anfray for his help with a Latin translation, and to an anonymous referee for pointing out various sources and exegetical aspects that were missing in the first version of this work. Any errors, whether historical or logical, are my own. This research was partly supported by the programs FRONTCOG (ANR-17-EURE-0017) and AMBISENSE (ANR-19-CE28-0019-01).}

\newpage

\bibliographystyle{apalike}
\bibliography{egre_biblio}

\newcommand{\SortNoop}[1]{}
\begin{thebibliography}{}

\bibitem[Ashworth, 1972a]{ashworth1972strict}
Ashworth, E.~J. (1972a).
\newblock Strict and material implication in the early {S}ixteenth {C}entury.
\newblock {\em Notre Dame Journal of Formal Logic}, 13(4):556--560.

\bibitem[Ashworth, 1972b]{ashworth1972treatment}
Ashworth, E.~J. (1972b).
\newblock The treatment of semantic paradoxes from 1400 to 1700.
\newblock {\em Notre Dame Journal of Formal Logic}, 13(1):34--52.

\bibitem[Ashworth, 1976]{ashworth1976bridge}
Ashworth, E.~J. (1976).
\newblock Will {S}ocrates cross the bridge? {A} problem in medieval logic.
\newblock {\em Franciscan Studies}, 36(1):75--84.

\bibitem[Barrio et~al., 2015]{barrio2015logics}
Barrio, E., Rosenblatt, L., and Tajer, D. (2015).
\newblock The logics of strict-tolerant logic.
\newblock {\em Journal of Philosophical Logic}, 44(5):551--571.

\bibitem[Beall, 2009]{beall2009spandrels}
Beall, J. (2009).
\newblock {\em Spandrels of truth}.
\newblock OUP Oxford.

\bibitem[Beaver, 2001]{beaver2001}
Beaver, D. (2001).
\newblock {\em Presupposition and assertion in dynamic semantics}, volume~29.
\newblock CSLI publications Stanford.

\bibitem[Biard, 1993]{biard1993sophismata}
Biard, J., editor (1993).
\newblock {\em Sophismes, de Jean Buridan}.
\newblock Jacques Vrin: Paris.
\newblock French edition and translation, based on a first translation by
  Fabienne Pironnet.

\bibitem[Bochvar, 1937]{bochvar1937}
Bochvar, D.~A. (1937).
\newblock [{O}n a three-valued calculus and its applications to the paradoxes
  of the classical extended functional calculus].
\newblock {\em Mathematicheskii sbornik}, 4(46):287--308.
\newblock English translation by M. Bergmann in \textit{History and Philosophy
  of Logic} 2 (1981), 87--112.

\bibitem[Bonevac, 1990]{bonevac1990paradoxes}
Bonevac, D. (1990).
\newblock Paradoxes of fulfillment.
\newblock {\em Journal of philosophical logic}, pages 229--252.

\bibitem[Burge, 1978]{burge1978buridan}
Burge, T. (1978).
\newblock Buridan and epistemic paradox.
\newblock {\em Philosophical Studies}, 34(1):21--35.

\bibitem[Cervantes, 1615]{quixote}
Cervantes, M. (1615).
\newblock {\em Don Quixote}.
\newblock Gutenberg project, 2019 edition.
\newblock English translation by John {O}rmsby.

\bibitem[Church, 1956]{church1996introduction}
Church, A. (1956).
\newblock {\em Introduction to mathematical logic}, volume~13.
\newblock Princeton University Press.
\newblock 1996 reissue.

\bibitem[Cobreros, 2016]{cobreros2016supervaluationism}
Cobreros, P. (2016).
\newblock Supervaluationism and the timeless solution to the foreknowledge
  problem.
\newblock {\em Scientia et Fides}, 4(1):61--75.

\bibitem[Cobreros et~al., 2012]{cobreros2012tcs}
Cobreros, P., Egr{\'{e}}, P., Ripley, D., and van Rooij, R. (2012).
\newblock Tolerant, classical, strict.
\newblock {\em The Journal of Philosophical Logic}, 41(2):347--385.

\bibitem[Cobreros et~al., 2013]{cobreros2013reaching}
Cobreros, P., {\'E}gr{\'e}, P., Ripley, D., and van Rooij, R. (2013).
\newblock Reaching transparent truth.
\newblock {\em Mind}, 122(488):841--866.

\bibitem[Cobreros et~al., 2015]{cobreros2015vagueness}
Cobreros, P., Egr\'e, P., Ripley, D., and van Rooij, R. (2015).
\newblock Vagueness, truth and permissive consequence.
\newblock In D.~Achouriotti, H.~Galinon, J.~M., editor, {\em Unifying the
  Philosophy of Truth}, pages 409--430. Springer.

\bibitem[Cobreros et~al., 2020]{cobreros2020inferences}
Cobreros, P., Egr{\'e}, P., Ripley, D., and van Rooij, R. (2020).
\newblock Inferences and metainferences in {S}{T}.
\newblock {\em Journal of Philosophical Logic}, 49(6):1057--1077.

\bibitem[Cooper, 1968]{cooper1968propositional}
Cooper, W. (1968).
\newblock The propositional logic of ordinary discourse.
\newblock {\em Inquiry}, 11(1-4):295--320.

\bibitem[Dicher and Paoli, 2019]{dicher2019st}
Dicher, B. and Paoli, F. (2019).
\newblock {S}{T}, {L}{P} and tolerant metainferences.
\newblock In Baskent, C. and Ferguson, T., editors, {\em Graham Priest on
  dialetheism and paraconsistency}, pages 383--407. Springer.

\bibitem[Dutilh~Novaes, 2009]{novaes2009lessons}
Dutilh~Novaes, C. (2009).
\newblock Lessons on sentential meaning from mediaeval solutions to the {L}iar
  paradox.
\newblock {\em The Philosophical Quarterly}, 59(237):682--704.

\bibitem[Dutilh~Novaes and Read, 2008]{dutilh2008insolubilia}
Dutilh~Novaes, C. and Read, S. (2008).
\newblock Insolubilia and the fallacy \textit{secundum quid et simpliciter}.
\newblock {\em Vivarium}, 46(2):175--191.

\bibitem[\'Egr\'e, 2019]{egre2019respects}
\'Egr\'e, P. (2019).
\newblock Respects for contradictions.
\newblock In Baskent, C. and Ferguson, T., editors, {\em Graham {P}riest on
  Dialetheism and Paraconsistency}, pages 39--57. Springer.

\bibitem[\'Egr\'e, 2021]{egre2021half}
\'Egr\'e, P. (2021).
\newblock Half-truths and the {L}iar.
\newblock In Nicolai, C. and Stern, J., editors, {\em Modes of Truth}, pages
  18--40. Routledge: Taylor and Francis.

\bibitem[\'Egr\'e et~al., 2021]{egre2021finettian}
\'Egr\'e, P., Rossi, L., and Sprenger, J. (2021).
\newblock De {F}inettian logics of indicative conditionals. {P}art {I}:
  Trivalent semantics and validity.
\newblock {\em Journal of Philosophical Logic}, 50:187--213.

\bibitem[Field, 2008]{field2008saving}
Field, H. (2008).
\newblock {\em Saving truth from paradox}.
\newblock OUP Oxford.

\bibitem[Glanzberg, 2004]{glanzberg2004contextual}
Glanzberg, M. (2004).
\newblock A contextual-hierarchical approach to truth and the {L}iar paradox.
\newblock {\em Journal of Philosophical Logic}, 33(1):27--88.

\bibitem[Herzberger, 1973]{herzberger1973dimensions}
Herzberger, H.~G. (1973).
\newblock Dimensions of truth.
\newblock {\em Journal of Philosophical Logic}, 2(4):535--556.

\bibitem[Hughes, 1982]{hughes1982buridan}
Hughes, G.~E. (1982).
\newblock {\em John {B}uridan on Self-Reference: Chapter Eight of {B}uridan's
  {S}ophismata, with a Translation, an Introduction, and a Philosophical
  Commentary}.
\newblock Cambridge University Press.

\bibitem[Hyde, 1997]{hyde1997heaps}
Hyde, D. (1997).
\newblock From heaps and gaps to heaps of gluts.
\newblock {\em Mind}, 106(424):641--660.

\bibitem[Jacquette, 1991]{jacquette1991buridan}
Jacquette, D. (1991).
\newblock Buridan's {B}ridge.
\newblock {\em Philosophy}, 66(258):455--471.

\bibitem[Jones, 1986]{jones1986liar}
Jones, J.~R. (1986).
\newblock The {L}iar {P}aradox in {D}on {Q}uixote {II}, 51.
\newblock {\em Hispanic Review}, 54(2):183--193.

\bibitem[Klima, 2001]{klimaburidan}
Klima, G., editor (2001).
\newblock {\em Summulae de dialectica, by Jean Buridan}.
\newblock Yale University Press.

\bibitem[Klima, 2008]{klima2008logic}
Klima, G. (2008).
\newblock Logic without truth.
\newblock In Rahman, S., Tulenheimo, T., and Genot, E., editors, {\em Unity,
  Truth and the Liar: The Modern Relevance of Medieval Solutions to the Liar
  Paradox}, pages 87--112. Springer.

\bibitem[Kretzmann and Kretzmann, 1990]{kilvingtonsophismata}
Kretzmann, N. and Kretzmann, B.~E., editors (1990).
\newblock {\em Sophismata, by Richard Kilvington}.
\newblock Cambridge University Press.

\bibitem[McGee, 1990]{mcgee1990truth}
McGee, V. (1990).
\newblock {\em Truth, vagueness, and paradox: An essay on the logic of truth}.
\newblock Hackett Publishing.

\bibitem[{\O}hrstr{\o}m and Hasle, 2007]{ohrstrom2007temporal}
{\O}hrstr{\o}m, P. and Hasle, P. (2007).
\newblock {\em Temporal logic: from ancient ideas to artificial intelligence},
  volume~57.
\newblock Springer Science \& Business Media.

\bibitem[Perini-Santos, 2011]{perini2011buridan}
Perini-Santos, E. (2011).
\newblock John {B}uridan's theory of truth and the paradox of the {L}iar.
\newblock {\em Vivarium}, 49(1-3):184--213.

\bibitem[Pironnet, 2004]{buridanpironnet}
Pironnet, F., editor (2004).
\newblock {\em Summulae de practica sophismatum, by Jean Buridan}.
\newblock Aristarium 10-9: Brepols.

\bibitem[Pironnet, 2008]{pironnet2008william}
Pironnet, F. (2008).
\newblock William {H}eytesbury and the treatment of \textit{{I}nsolubilia} in
  fourteenth-century {E}ngland followed by a critical edition of three
  anonymous treatises \textit{{D}e insolubilibus} inspired by {H}eytesbury.
\newblock In Rahman, S., Tulenheimo, T., and Genot, E., editors, {\em Unity,
  Truth and the Liar: The Modern Relevance of Medieval Solutions to the Liar
  Paradox}, pages 255--333. Springer.

\bibitem[Priest, 1979]{priest1979logic}
Priest, G. (1979).
\newblock The logic of paradox.
\newblock {\em Journal of Philosophical logic}, 8(1):219--241.

\bibitem[Priest, 2006]{priest2006incontradiction}
Priest, G. (2006).
\newblock {\em In Contradiction}.
\newblock Oxford University Press.
\newblock Second edition.

\bibitem[Priest, 2019]{priest2019respectfully}
Priest, G. (2019).
\newblock Respectfully yours.
\newblock In Baskent, C. and Ferguson, T., editors, {\em Graham {P}riest on
  Dialetheism and Paraconsistency}. Springer.

\bibitem[Priest, 2021]{priest2021substructural}
Priest, G. (2021).
\newblock Substructural solutions to the semantic paradoxes: Dialetheism in
  sheepÕs clothing?
\newblock Manuscript.

\bibitem[Prior, 1962]{prior1962problems}
Prior, A. (1962).
\newblock Some problems of self-reference in {J}ohn {B}uridan.
\newblock pages 281--296.
\newblock Dawes Hicks Lecture on Philosophy.

\bibitem[Quine, 1953]{quine1953so}
Quine, W. V.~O. (1953).
\newblock On a so-called paradox.
\newblock {\em Mind}, 62(245):65--67.

\bibitem[Read, 1995]{read1995thinking}
Read, S. (1995).
\newblock {\em Thinking about Logic}.
\newblock Oxford University Press.

\bibitem[Read, 2002]{read2002liar}
Read, S. (2002).
\newblock The {L}iar paradox from {J}ohn {B}uridan back to {T}homas
  {B}radwardine.
\newblock {\em Vivarium}, 40(2):189--218.

\bibitem[Read, 2010]{bradwardineinsol}
Read, S., editor (2010).
\newblock {\em Insolubilia, by Thomas Bradwardine}.
\newblock Dallas Medieval Texts and Translations 10: Peeters.
\newblock Edition, translation and introduction by Stephen Read.

\bibitem[Ripley, 2012]{ripley2012conservatively}
Ripley, D. (2012).
\newblock Conservatively extending classical logic with transparent truth.
\newblock {\em The Review of Symbolic Logic}, 5(2):354--378.

\bibitem[Ripley, 2013]{ripley2013paradoxes}
Ripley, D. (2013).
\newblock Paradoxes and failures of cut.
\newblock {\em Australasian Journal of Philosophy}, 91(1):139--164.

\bibitem[Rossi, 2019]{rossi2019unified}
Rossi, L. (2019).
\newblock A unified theory of truth and paradox.
\newblock {\em The Review of Symbolic Logic}, 12(2):209--254.

\bibitem[Scott, 1966]{buridanscott}
Scott, T.~K., editor (1966).
\newblock {\em Sophisms on Meaning and Truth, by Jean Buridan}.
\newblock Appleton-Century-Crofts: New York.

\bibitem[Spade and Wilson, 1986]{wyclifsumma}
Spade, P.~V. and Wilson, G.~A., editors (1986).
\newblock {\em Johannis Wyclif Summa Insolubilium}.
\newblock Medieval and Renaissance Texts and Studies: Binghampton.

\bibitem[Spector, 2016]{spector2016multivalent}
Spector, B. (2016).
\newblock Multivalent semantics for vagueness and presupposition.
\newblock {\em Topoi}, 35(1):45--55.

\bibitem[Thomason, 1970]{thomason1970indeterminist}
Thomason, R.~H. (1970).
\newblock Indeterminist time and truth-value gaps.
\newblock {\em Theoria}, 36(3):264--281.

\bibitem[Ulatowski, 2003]{ulatowski2003buridan}
Ulatowski, J.~W. (2003).
\newblock A conscientious resolution of the action paradox on {B}uridan's
  bridge.
\newblock {\em Southwest Philosophical Studies}, pages 85--93.

\bibitem[Van~Fraassen, 1968]{fraassen1968presupposition}
Van~Fraassen, B. (1968).
\newblock Presupposition, implication, and self-reference.
\newblock {\em The Journal of Philosophy}, 65(5):136--152.

\bibitem[Walker, 1847]{murray1847compendium}
Walker, J., editor (1847).
\newblock {\em Murray's Compendium of Logic}.
\newblock Longman, Brown, Green, Longmans: London.

\bibitem[Zehr, 2014]{Zehr2014:PHD}
Zehr, J. (2014).
\newblock {\em Vagueness, Presupposition and Truth-Value Judgments}.
\newblock PHD Thesis, ENS, PSL {U}niversity, Paris.

\end{thebibliography}

\end{document}